\documentclass[a4paper,12pt]{article}

\usepackage{amssymb,amsmath,amsthm,amscd,latexsym}

\newcommand{\Om} {\Omega}
\newcommand{\pa} {\partial}
\newcommand{\al} {\alpha}

\newcommand{\De} {\Delta}

\newcommand{\noi} {\noindent}

\newcommand{\ep} {\epsilon}
\newcommand{\na} {\nabla}
\newcommand{\be} {\begin{equation}}
\newcommand{\ee} {\end{equation}}
\newcommand{\la} {\lambda}
\newcommand{\de} {\delta}

\def\lan{\langle}
\def\ran{\rangle}

\def\Ga{\Gamma}
\def\Xint#1{\mathchoice
   {\XXint\displaystyle\textstyle{#1}}%
   {\XXint\textstyle\scriptstyle{#1}}%
   {\XXint\scriptstyle\scriptscriptstyle{#1}}%
   {\XXint\scriptscriptstyle\scriptscriptstyle{#1}}%
   \!\int}
\def\XXint#1#2#3{{\setbox0=\hbox{$#1{#2#3}{\int}$}
     \vcenter{\hbox{$#2#3$}}\kern-.5\wd0}}

\def\aint{\Xint\diagup}

\begin{document}
\title{ A note on Alexsandrov type theorem for $k$-convex functions}

\author{{\Large Nirmalendu Chaudhuri\,and\, Neil S. Trudinger} \\\\ Centre for Mathematics 
and its Applications\\ 
Australian National University\\
Canberra, ACT 0200\\
Australia\\
chaudhur@maths.anu.edu.au\\
 neil.trudinger@maths.anu.edu.au
}
\maketitle
\begin{center}{Abstract}
\medskip

{\footnotesize\rm \begin{tabular}{p{14.5cm}}
In this note we show that $k$-convex functions on $\Bbb R^n$ are twice 
differentiable almost everywhere for every positive integer $k>n/2$. 
This generalizes the classical Alexsandrov's theorem for convex functions. 
\end{tabular} } 
\end{center}

\section{Introduction} 

A classical result of Alexsandrov \cite{A} asserts that convex functions in $\Bbb R^n$ 
are twice differentiable $a.e.$, (see also \cite{EG}, \cite{Kr} for more modern 
treatments). It is well known that Sobolev functions $u\in W^{2,p}$, for $p>n/2$ are 
twice differentiable a.e.. The following weaker notion of convexity known as $k$-convexity 
was introduced by Trudinger and Wang \cite{TW1, TW2}. Let $\Om\subset\Bbb R^n$ be an open set and $C^2(\Om)$ be the 
 class of continuously twice differentiable functions on $\Om$. For $k=1,\,2,\,\dots,n$ and a function 
 $u\in C^2(\Om)$, the $k$-Hessian operator, $F_k$, is defined by      
 
 \begin{equation}
 \label{a1}
 F_k[u]:=\,S_k(\la(\na^2u))\,,
 \end{equation}  
 where $\na^2u=(\pa_{ij}u)$ denotes the Hessian matrix of the second derivatives of $u$, 
 $\la(A)\,=\,(\la_1,\la_2,\cdots,\la_n)$ the vector of eigenvalues of an $n\times n$ 
 matrix $A\in\Bbb M^{n \times n}$ and $S_k(\la)$ is the $k$-th elementary symmetric 
 function on $\Bbb R^n$, 
 given by 

 \begin{equation}
 \label{a2}
 S_k(\la):=\sum_{i_1<\cdots<i_k}\la_{i_1}\cdots \la_{i_k}\,.
 \end{equation}

Alternatively we may write 
\begin{equation}
\label{a3}
F_k[u]=\,[\na^2u]_k\,
\end{equation}
 where $[A]_k$ denotes the sum of $k\times k$ principal minors of an $n\times n$ matrix 
$A$, which may also be called the $k$-trace of $A$. The study of $k$-Hessian operators was 
initiated by Caffarelli, Nirenberg and Spruck \cite{CNS} and Ivochkina \cite{I} with
 further developed by Trudinger and Wang \cite{T, TW1, TW2, TW4, TW3}. 
 
A function $u\in C^2(\Om)$ is called $k$-{\it convex} in $\Om$ if 
$F_j[u]\,\geq\,0$ in $\Om$ for $j=1,\,2\,\dots, k$; that is, the eigenvalues $\la(\na^2u)$ of the 
Hessian $\na^2u$ of $u$ lie in the closed convex cone given by 
\begin{equation}
\label{a4}
\Gamma_k\,:=\,\{\la\in\Bbb R^n\,\,: \,\,S_j(\la)\,\geq\,0,\,\,j=1,\,2,\,\dots, k\}\,.
\end{equation}  
(see \cite{CNS} and \cite{TW2} for the basic properties of $\Gamma_k$.) We notice that 
$F_1[u]\,=\,\De u$, is the Laplacian operator and $1$-convex functions are subharmonic. When $k=n$,
$F_n[u]\,=\,{\rm det}(\na^2u)$, the Monge-Amp\'ere operator and $n$-convex functions are convex. 
To extend  the definition of $k$-convexity for non-smooth functions we adopt a {\it viscosity} 
definition as in \cite{TW2}. {\it An upper semi-continuous function 
$u\,:\,\Om\to [-\infty,\infty)$ ($u\not\equiv -\infty$ on any 
connected component of $\Om$) is called $k$-convex if $F_j[q]\,\geq\,0$, in $\Om$ for $j=1,2,\dots,k$, 
for every quadratic polynomial $q$ for which the difference $u-q$ has a finite local maximum in $\Om$}. 
Henceforth, we shall denote the class of $k$-convex functions in $\Om$ by $\Phi^k(\Om)$.  
When $k=1$ the above definition is equivalent to the usual 
definition of subharmonic function, see, for example (Section 3.2, \cite{H}) or 
(Section 2.4, \cite{ MK}). Thus $\Phi^1(\Om)$ is the class of subharmonic functions in $\Om$. 
We notice that $\Phi^k(\Om)\subset \Phi^1(\Om)\subset L^1_{\rm loc}(\Om)$ for $k=1,2,\dots,n$, 
and a function $u\in \Phi^n(\Om)$ if and only if it is convex on each component of $\Om$.
Among other results Trudinger and Wang \cite{TW2} (Lemma 2.2) proved that $u\in\Phi^k(\Om)$ if and only if 
\begin{equation}
\label{a9} 
\int_{\Om}u(x)\,\left(\sum_{i,j}^{n}a^{ij}\pa_{ij}\phi(x)\right)\,dx\geq\,0
\end{equation}
for all smooth compactly supported functions $\phi\,\geq\,0$, and for all constant $n\times n $ 
symmetric matrices $A=(a^{ij})$ with eigenvalues $\la(A)\in \Ga_k^*$, where $\Ga_k^*$ is dual 
cone defined by  
\begin{equation}
\label{a7} 
\Ga_k^*:\,=\,\{\la\in\Bbb R^n\,\,:\,\,\lan \la,\mu\ran\,\geq 0
\,\,\,\,{\rm for}\,\,\,{\rm all}\,\,\,\mu\in \Ga_k\}\,. 
\end{equation}
In this note we prove the following Alexsandrov type theorem for $k$-convex functions. 
\medskip 

\noi
{\bf Theorem 1.1.} \, {\it Let $k>n/2$, $n\geq 2$ and $u\,:\,\Bbb R^n\to [-\infty,\infty)$ ($u\not\equiv-\infty$ 
on any component of $\Bbb R^n$), be a $k$-convex function. Then $u$ is twice differentiable almost 
everywhere. More precisely, we have the Taylor's series expansion for ${\cal L}^n$ $x$ a.e., 
\begin{equation}
\label{a10}
\left|u(y)-u(x)-\lan \na u(x)\,, \,y-x\ran -\frac{1}{2}\lan \na^2u(x)(y-x)\,,\,y-x\ran\right|
=o(|y-x|^2)\,,
\end{equation}
as $y\to x$.}
\medskip

In Section 3 (see, Theorem 3.2.), we also prove that the absolute continuous part of the $k$-Hessian measure 
(see, \cite{TW1, TW2}) $\mu_k[u]$, associated to a $k$-convex function for $k>n/2$ is represented by 
$F_k[u]$. For the Monge-Amp\'ere measure $\mu[u]$ associated to a convex function $u$, such 
result is obtained in \cite{TW5}. 
\medskip

To conclude this introduction we note that it is equivalent to assume only $F_k[q]\geq 0$, 
in the definition of $k$-convexity \cite{TW2}. Moreover $\Gamma_k$ may also be characterized as the 
closure of the positivity set of $S_k$ containing the positive cone $\Gamma_n$, \cite{CNS}.

\section{Notations and preliminary results}\setcounter{equation}{0} 

Throughout the text we use following standard notations. $|\cdot|$ and 
$\lan\cdot,\cdot\ran$ 
will stand for Euclidean norm and inner product in $\Bbb R^n$, and   
$B(x,r)$ will denote the open ball in $\Bbb R^n$ of radius $r$ centered at $x$. For measurable 
 $E\subset\Bbb R^n$, ${\cal L}^n$(E) will denote its Lebesgue measure. For a smooth function 
$u$, the gradient and Hessian of $u$ are denoted by ${\displaystyle\na u=(\pa_1 u,\cdots,\pa_n u)}$ and 
${\displaystyle \na^2u=(\pa_{ij} u)_{1\leq i,j\leq n}}$ respectively. For a locally integrable 
function $f$, the distributional gradient and Hessian are denoted 
by $Df=(D_1f,\cdots, D_nf)$ and ${\displaystyle D^2u=(D_{ij} u)_{1\leq i,j\leq n}}$ respectively. 
\medskip

For the convenience of the readers, we cite the following H\"older and gradient 
estimates for $k$-convex functions, and the weak continuity result for $k$-Hessian measures, 
\cite{TW1, TW2}.  
\medskip

\noindent
{\bf Theorem 2.1.} (Theorem 2.7, \cite{TW2})\, {\it For $k\,>\,n/2$, $\Phi^k(\Om)\subset 
C^{0,\al}_{\rm loc}(\Om)$ with $\al:=2-n/k$ and for any $\Om^{\prime}\subset\subset\Om$, 
$u\in \Phi^k(\Om)$, there exists $C>0$, depending only on $n$ and $k$ such that 
\begin{equation} 
\label{a6}
\sup_{{\displaystyle\stackrel {x,y\in \Om^{\prime}}{x\neq y}}}
d_{x,y}^{n+\al}\frac{|u(x)-u(y)|}{|x-y|^{\al}}\,\leq\,C\,\int_{\Om^{\prime}} |u|\,\,,
\end{equation}
where $d_x\,:=\,{\rm dist}(x\,,\,\pa\Om^{\prime})$ and $d_{x,y}\,:=\,\min\{d_x\,,\,d_y\}$.}
\medskip

\noindent 
{\bf Theorem 2.2.} (Theorem 4.1, \cite{TW2})\, {\it For $k=1,\dots, n$,\, and\, 
${\displaystyle 0\,<\,q\,<\,\frac{nk}{n-k}}$, the space of $k$-convex functions 
$\Phi^k(\Om)$ lie in the local Sobolev space $W^{1,q}_{\rm loc}(\Om)$. Moreover, for any $\Om^{\prime}
\subset\subset\Om^{\prime\prime}\subset\subset\Om$ and $u\in\Phi^k(\Om)$ there exists $C>0$, 
depending on $n$, $k$, $q$, $\Om^{\prime}$ and $\Om^{\prime\prime}$,  
 such that
 \begin{equation}
 \label{a12}
 \left(\int_{\Om^{\prime}}|Du|^q\right)^{1/q}\,\leq\, C\int_{\Om^{\prime\prime}} |u|\,.
\end{equation}
}
\medskip

\noi 
 {\bf Theorem 2.3.} [Theorem 1.1, \cite{TW2}]\, {\it For any $u\in \Phi^k(\Om)$, there exists a 
Borel measure $\mu_k[u]$ in $\Om$ such that
\medskip

\noindent
{\rm (i)}\, $\mu_k[u](V)=\int_V F_k[u](x)\,dx$\, for any Borel set $V\subset \Om$, if $u\in C^2(\Om)$ and 
\medskip

\noindent
{\rm (ii)} if $(u_m)_{m\geq 1}$ is a sequence in $\Phi^k(\Om)$ converges in $L^1_{\rm loc}(\Om)$ to 
a function $u\in \Phi^k(\Om)$, the sequence of Borel measures $(\mu_k[u_m])_{m\geq 1}$ converges weakly 
to $\mu_k[u]$.}
 \medskip
 
 Let us recall the definition of the dual cones, \cite{T2}
$$
\Ga_k^*:\,=\,\{\la\in\Bbb R^n\,\,:\,\,\lan \la,\mu\ran\,\geq 0
\,\,\,\,{\rm for}\,\,\,{\rm all}\,\,\,\mu\in \Ga_k\}\,, 
$$
which are also closed convex cones in $\Bbb R^n$. We notice that $\Ga_j^*\subset \Ga_k^*$ for 
$j\,\leq\, k$ with $\Ga_n^*=\Ga_n=\{\la\in\Bbb R^n\,\,: \,\,\la_i
\geq\,0,\,\,\,\,j=1,\,2,\,\dots, n\}$, $\Ga_1^*$ is the ray given by 
$$
\Ga
_1^*\,=\,\{\,t(1,\,\cdots,\,1)\,\,:\,\, t\,\geq\,0\}\,,
$$ 
and $\Ga_2^*$ has the following interesting characterization, 
\begin{equation}
\label{a8}
\Ga_2^*\,=\,\left\{\la\in \Ga_n\,\,: \,\,|\la|^2\,\leq\,\frac{1}{n-1}\left(\sum_{i=1}^n\la_i
\right)^2\,\right\}\,.
\end{equation}
We use this explicit representation of $\Ga_2^*$ to establish that the distributional 
derivatives $D_{ij}u$ of the $k$-convex function $u$ are signed Borel measures
 for $k\geq 2$, (see also \cite{TW2}). 
\medskip 

\noindent
{\bf Theorem 2.4.}\, {\it Let $2\leq k\leq n$ and $u\,:\,\Bbb R^n\to [-\infty,\infty)$, be a 
$k$-convex function. Then there exist signed Borel measures $\mu^{ij}=\mu^{ji}$ such that 
\begin{equation}
\label{a11}
\int_{\Bbb R^n}u(x)\,\pa_{ij}\phi(x)\,dx=\int_{\Bbb R^n}\phi(x)\, d\mu^{ij}(x)\,\,
\,\,\,\,{\rm for}\,\,\,i,\,j=1,2,\,\dots, n\,,
\end{equation}
for all $\phi\in C^{\infty}_c(\Bbb R^n)$.}
\medskip

 {\it Proof.}\, Let $k\geq 2$ and $u\in\Phi^k(\Bbb R^n)$. Since $\Phi^k(\Bbb R^n)
 \subset \Phi^2(\Bbb R^n)$ for $k\geq 2$, it is enough to prove the theorem for $k=2$. Let $u$ 
 be a $2$-convex function in $\Bbb R^n$. For $A\in \Bbb S^{n\times n}$, the space of 
 $n\times n$ symmetric matrices,  define the distribution 
 $T_A\,:\,C^2_c(\Bbb R^n)\to \Bbb R$, by 
 $$
 T_A(\phi)\,:=\int_{\Bbb R^n}u(x)\sum_{i,j}^{n} a^{ij}\pa_{ij}\phi(x)\,dx
 $$
By (\ref{a9}), $T_A(\phi)\,\geq 0$ for $A\in \Bbb S^{n\times n}$ with eigenvalues 
$\la(A)\in \Ga_2^*$, and $\phi\geq 0$. Therefore, by Riesz representation 
(see, for example Theorem 2.14 in \cite{Ru} or Theorem 1, Section 1.8 in \cite{EG}),
 there exist a Borel measure $\mu^A$ in $\Bbb R^n$, such that
\begin{equation}
\label{a13}
 T_A(\phi)=\int_{\Bbb R^n}\phi\,\sum_{i,j}^{n} a^{ij}D_{ij}u\,dx
 =\int_{\Bbb R^n}\phi\,d\mu^A\,,
\end{equation}
for all $\phi\in C^2_c(\Bbb R^n)$ and all $n\times n$ symmetric matrices $A$ with $\la(A)\in\Ga_2^*$. 
In order to prove the second order distributional derivatives $D_{ij}u$ of $u$ to be signed 
Borel measures, we need to make special choices for the matrix $A$. By taking $A=I_n$, the identity 
matrix, $\la(A)\in\Ga_1^*\subset\Ga_2^*$, we obtain a Borel measure $\mu^{I_n}$ such that   
\begin{equation}
\label{a14}
\int_{\Bbb R^n}\phi\,\sum_{i=1}^{n} D_{ii}u\,dx
 =\int_{\Bbb R^n}\phi\,d\mu^{I_n}\,,
\end{equation}
for all $\phi\in C^2_c(\Bbb R^n)$. Therefore, the trace of the distributional Hessian $D^2u$, 
is a Borel measure. For each $i=1,\dots, n$, let 
$A_i$ be the diagonal matrix with all entries $1$ but the $i$-th diagonal entry being $0$. 
Then by the characterization of $\Ga_2^*$ in $(\ref{a8})$, it follows that $\la(A_i)\in \Ga_2^*$. 
Hence there exist a Borel measure $\mu^i$ in $\Bbb R^n$ such that 
 \begin{equation}
\label{a15}
 \int_{\Bbb R^n}\phi\,\sum_{j\neq i}^{n} D_{jj}u\,dx
 =\int_{\Bbb R^n}\phi\,d\mu^{i}\,,
\end{equation}
for all $\phi\in C^2_c(\Bbb R^n)$. From $(\ref{a14})$ and $(\ref{a15})$ it follows that, the 
diagonal entries $D_{ii}u=\mu^{I_n}-\mu^{i}:=\mu^{ii}$ are signed Borel 
measure and 
\begin{equation}
\label{a16}
\int_{\Bbb R^n}u\,\pa_{ij}\phi\,dx =\int_{\Bbb R^n}\phi\,d\mu^{ii}\,,
\end{equation}
for all $\phi\in C^2_c(\Bbb R^n)$. Let $\{e_1,\dots, e_n\}$ be the standard 
orthonormal basis in $\Bbb R^n$ and for $a,\,b\in\Bbb R^n$, $a\otimes b:=(a^i b^j)$, denotes 
the $n\times n$ rank-one matrix. For $0<t<1$ and $i\neq j$, let us 
define $A_{ij}:=I_n +t [\,e_i\otimes e_j+e_j\otimes e_i]$. By a straight forward calculation, it is easy to see that 
the vector of eigenvalues $\la(A_{ij})=(1-t,1+t,1,\cdots,1)\in \Ga_2^*$, for 
${\displaystyle 0<t<\left(n/2(n-1)\right)^{1/2}}$\,. Note that for this choice of   
$A_{ij}$  
$$
\sum_{k,l=1}^{n} a^{kl}\pa_{kl}\phi=\sum_{k=1}^{n}\pa_{kk}\phi + 2t\, 
\pa_{ij}\phi\,.
$$
Thus for $i\neq j$, $(\ref{a13})$ and $(\ref{a14})$ yields
\begin{align}
\label{a17}
\int_{\Bbb R^n}u\,\pa_{ij}\phi\,dx
&=\frac{1}{2t}\left[\int_{\Bbb R^n}u\,\sum_{k,l=1}^{n} a^{kl}\pa_{kl}\phi
\,dx -\int_{\Bbb R^n}u\,\sum_{k=1}^{n}\pa_{kk}\phi\,dx\right]\nonumber\\
&=\frac{1}{2t}\left[\int_{\Bbb R^n}\phi\,d\mu^{A_{ij}} -\int_{\Bbb R^n}\phi\,d\mu^{I_n}
\right]\nonumber\\
&=\int_{\Bbb R^n}\phi\,d\mu^{ij}\,,
\end{align}
where 
$$
\mu^{ij}:=\frac{1}{2t}\left(\mu^{A_{ij}}-\mu^{I_n}\right)
=\frac{1}{2t}\left(\mu^{A_{ij}}-\sum_{k=1}^{n}\mu^{kk}\right)\,.
$$
Therefore $D_{ij}u=\mu^{ij}$, are signed Borel measures and satisfies the identity 
$(\ref{a11})$.\qed
 \medskip

A function $f\in L^1_{\rm loc}(\Bbb R^n)$ is said to have 
{\it locally bounded variation} in $\Bbb R^n$ if for each bounded open subset 
$\Om^{\prime}$ of $\Bbb R^n$, 
$$
\sup\left\{\int_{\Om^{\prime}}f\,{\rm div}\,\phi\,dx\,: \,\phi \in C^1_{c}(\Om^{\prime};\Bbb R^n),\,\,
|\phi(x)|\leq 1\,\,{\rm for}\,\,{\rm all}\,\,x\in\Om^{\prime}\right\}\,<\,\infty\,.
$$
 We use the notation $BV_{\rm loc}(\Bbb R^n)$, to denote the space of such functions. 
 For the theory of {\it functions of bounded variation} readers are referred to \cite{Gi, Zi, EG}. 
 \bigskip
 
 \noindent
 {\bf Theorem 2.5.}\, {\it Let $n\geq 2$, $k> n/2$ and $u\,:\,\Bbb R^n\to [-\infty,\infty)$, be a 
$k$-convex function. Then $u$ is differentiable a.e. ${\cal L}^n$ and ${\displaystyle
\frac{\pa u}{\pa x_i}\in BV_{\rm loc}(\Bbb R^n)}$, for all $i=1,\dots, n$.}
\bigskip

{\it Proof.} Observe that for $k>n/2$, we can take ${\displaystyle n\,<\,q\,<\,\frac{nk}{n-k}}$ and by the 
gradient estimate (\ref{a12}), we conclude that $k$-convex functions are differentiable 
${\cal L}^n$ a.e. $x$. Let $\Om^{\prime}\subset\subset\Bbb R^n$, 
$\phi=(\phi^1,\cdots,\phi^n)\in C^1_c(\Om^{\prime};\Bbb R^n)$ such that 
$|\phi(x)|\leq 1$ for $x\in\Om^{\prime}$. Then by integration by parts and 
the identity $(\ref{a11})$, we have for $i=1,\dots, n$, 
\begin{align}
\int_{\Om^{\prime}}\frac{\pa u}{\pa x_i}\,{\rm div}\,\phi\,dx
&=-\sum_{j=1}^{n}\int_{\Om^{\prime}}u
\frac{\pa^2 \phi^j}{\pa x_i\pa x_j}\,dx\nonumber\\
&=-\sum_{j=1}^{n}\int_{\Om^{\prime}}\phi^{j}\,d\mu^{ij}\nonumber\\
&\leq \sum_{j=1}^{n}|\mu^{ij}|(\Om^{\prime})<\infty\,,\nonumber
\end{align}
where $|\mu^{ij}|$ is the total variation of the Radon measure $\mu^{ij}$. 
This proves the theorem.\qed 

\section{Twice differentiability}\setcounter{equation}{0}

Let $u$ be a $k$-convex function, $k\geq 2$, then by the Theorem 2.4, we have 
 $D^2u=(\mu^{ij})_{i,j}$, where $\mu^{ij}$ are Radon measures. 
By Lebesgue's Decomposition Theorem, we may write 
$$
\mu^{ij}=\mu_{\rm ac}^{ij}+\mu^{ij}_{\rm s}\,\,\,\,\,\,\,\,{\rm for}\,\,\,i,j=1,\cdots,n\,,
$$
where $\mu_{\rm ac}^{ij}$ is absolutely continuous with respect to ${\cal L}^n$ and 
$\mu^{ij}_{\rm s}$ is supported on a set with Lebesgue measure zero.  
Let $u_{ij}$ be the density of the absolutely continuous part, i.e., 
$d\mu_{\rm ac}^{ij}=u_{ij}\, dx$, $u_{ij}\in L^1_{\rm loc}(\Bbb R^n)$. Set 
${\displaystyle u_{ij}:=\frac{\pa^2 u}{\pa x_i \pa x_j}}$,\,  
${\displaystyle \na^2 u:= \left(\frac{\pa^2 u}{\pa x_i \pa x_j}\right)_{i,j}=(u_{ij})_{i,j}\in 
L^1_{\rm loc}(\Bbb R^n;\Bbb R^{n\times n})}$ and 
$[D^2u]_s:=(\mu^{ij}_{\rm s})_{i,j}$. Thus the vector valued Radon measure $D^2u$ 
can be decomposed as $D^2u=[D^2u]_{\rm ac} 
+[D^2u]_s$, where $d[D^2u]_{\rm ac}=\na^2u\,dx$. Now we are in a position to 
prove the theorem 1.1. To carry out the proof, we use a similar 
approach to Evans and Gariepy, see, Section 6.4, in \cite{EG}. 
\bigskip

\noindent
{\it Proof of Theorem 1.1.}\, 
 Let $n\geq 2$ and $u$ be a $k$-convex function on $\Bbb R^n$, 
$k>n/2$. Then by Theorem 2.4, and Theorem 2.5, we have for ${\cal L}^n$ a.e. $x$
\begin{equation}
\label{a18}
\lim_{r\to 0}\aint_{B(x,r)}|\na u(y)-\na u(x)|\,dy\,=\,0\,,
\end{equation}
\begin{equation}
\label{a19}
\lim_{r\to 0}\aint_{B(x,r)}|\na^2 u(y)-\na^2 u(x)|\,dy\,=\,0\,
\end{equation}
and 
\begin{equation}
\label{a20}
\lim_{r\to 0}\frac{|[D^2u]_s|(B(x,r))}{r^n}\,=\,0\,.
\end{equation}
where $\aint_{E} f\,dx$ we denote the mean value 
$\left({\cal L}^n(E)\right)^{-1}\int_E f\,dx$. Fix a point x for which (\ref{a18})-(\ref{a20}) holds.
Without loss generality we may assume $x=0$. Then following similar calculations as in the 
proof of Theorem 1, Section 6.4 in \cite{EG}, we obtain, 
\begin{equation}
\label{a26}
\aint_{B(r)}\left|u(y)-u(0) -\left\lan \na u(0),y\right\ran -\frac{1}{2}
\left\lan \na^2u(0)y,y\right\ran\right|\,dy=o(r^2)\,,
\end{equation}
as $r\to 0$. In order to establish 
\begin{equation}
\label{a27}
\sup_{B(r/2)}\left|u(y)-u(0) -\left\lan \na u(0),y\right\ran -\frac{1}{2}
\left\lan \na^2u(0)y,y\right\ran\right|=o(r^2)\,\,\,\,{\rm as}\,\,r\to 0\,,
\end{equation}
we need the following lemma.
\medskip 

\noi
{\bf Lemma 3.1.}\, {\it Let ${\displaystyle h(y):=u(y)-u(0) -\left\lan \na u(0),y\right\ran -\frac{1}{2}
\left\lan \na^2u(0)y,y\right\ran}$. Then there exists a constant $C>0$ depending only 
on $n$, $k$ and $|\na^2 u(0)|$, such that for any $0<r<1$
\begin{equation}
\label{a28}
\sup_{{\displaystyle\stackrel {y,z\in B(r)}{y\neq z}}}\frac{|h(y)-h(z)|}{|y-z|^{\al}}
\leq \frac{C}{r^{\al}}\aint_{B(2r)}|h(y)|\,dy +Cr^{2-\al}\,,
\end{equation}
where $\al:=(2-n/k)$}. 
\medskip

{\it Proof.}\, Let $\Lambda :=|\na^2u(0)|$ and define 
${\displaystyle g(y):=h(y)+\frac{\Lambda}{2} |y|^2}$. Since 
$\frac{\Lambda}{2} |y|^2-u(0) -\left\lan \na u(0),y\right\ran -\frac{1}{2}
\left\lan \na^2u(0)y,y\right\ran$ is convex and sum of two $k$-convex functions 
are $k$-convex (follows from (\ref{a4})), 
we conclude that $g$ is $k$-convex. Applying the H\"older 
estimate in (\ref{a6}) for $g$ with $\Om^{\prime}=B(2r)$, there exists $C:=C(n,k)>0$, such that
\begin{align}
\label{a29}
r^{n+\al}\sup_{{\displaystyle\stackrel {y,z\in B(r)}{y\neq z}}}\frac{|g(y)-g(z)|}{|y-z|^{\al}}
&={\rm dist}(B(r),\pa B(2r))^{n+\al}
\sup_{{\displaystyle\stackrel {y,z\in B(r)}{y\neq z}}}\frac{|g(y)-g(z)|}{|y-z|^{\al}}\nonumber\\
&\leq \sup_{{\displaystyle\stackrel {y,z\in B(2r)}{y\neq z}}}d_{y,z}^{n+\al}\frac{|g(y)-g(z)|}{|y-z|^{\al}}\nonumber\\
&\leq C\int_{B(2r)}|g(y)|\,dy\nonumber\\
&\leq C\int_{B(2r)}|h(y)|\,dy+ Cr^{n+2}\,,
\end{align} 
where $d_{y,z}:=\min\{{\rm dist}(y,\pa B(2r))\,,\,{\rm dist}(z,\pa B(2r))\}$. 
Therefore the estimate (\ref{a28}) for $h$ follows from the estimate 
(\ref{a29}) and the definition of $g$.\qed     
\medskip

{\it Proof of Theorem 1.1. (ctd.)} To prove (\ref{a27}), take $0\,<\,\ep,\,\delta\,<\,1$, such that $\delta^{1/n}\,\leq\, 1/2$. 
Then there exists $r_0$ depending on $\ep$ and $\delta$, sufficiently small, 
 such that, for $0\,<\,r\,<\,r_0$ 
\begin{align}
\label{a30}
{\cal L}^n\left\{z\in B(r)\,\,:\,\,|h(z)|\,\geq\,\ep r^2\right\}
&\leq \frac{1}{\ep r^2}\int_{B(r)}|h(z)|\,dz\nonumber\\
&=o(r^n)~~~~~~{\rm by}~~~(\ref{a26})\nonumber\\
&<\delta {\cal L}^n(B(r))
\end{align}
Set $\sigma:=\delta^{1/n}r$. Then for each $y\in B(r/2)$ there exists $z\in B(r)$ such that 
$$
|h(z)|\leq \ep r^2\,~~~~{\rm and}~~~~
|y-z|\leq\sigma.
$$
 Hence for each $y\in B(r/2)$, we obtain by 
(\ref{a26}) and (\ref{a28}),  
\begin{align}
|h(y)|&\leq |h(z)| +|h(y)-h(z)|\nonumber\\
&\leq \ep r^2 + C|y-z|^{\al}\left(\frac{1}{r^{\al}}\aint_{B(2r)}|h(y)|\,dy +r^{2-\al}\right)\nonumber\\
&\leq \ep r^2 + C\delta^{\al/n}r^{\al}\left(\frac{1}{r^{\al}}\aint_{B(2r)}|h(y)|\,dy 
+r^{2-\al}\right)\nonumber\\
&\leq \ep r^2 + C\delta^{\al/n}\left(\aint_{B(2r)}|h(y)|\,dy 
+r^{2}\right)\nonumber\\
&=r^2\left(\ep +C\delta^{\al/n}\right) +o(r^2)~~~~~{\rm as}~~~r\to 0\nonumber
\end{align}
By choosing $\delta$ such that, $ C\delta^{\al/n}=\ep$, we have for sufficiently small $\ep>0$ and 
$0<r<r_0$, 
$$
\sup_{B(r/2)}|h(y)|\leq 2\ep r^2 +o(r^2)\,.
$$
Hence 
$$
\sup_{B(r/2)}\left|u(y)-u(0) -\left\lan \na u(0),y\right\ran -\frac{1}{2}
\left\lan \na^2u(0)y,y\right\ran\right|\,dy=o(r^2)\,\,\,\,{\rm as}\,\,r\to 0\,.
$$
This proves (\ref{a10}) for $x=0$ and hence $u$ is twice differentiable at $x=0$. 
Therefore $u$ is twice differentiable at every $x$ and satisfies (\ref{a10}), for which 
(\ref{a18})-(\ref{a20}) holds. This proves the theorem.\qed 
\bigskip

Let $u$ be a $k$-convex function and $\mu_k[u]$ be the associated $k$-Hessian measure. 
Then $\mu_k[u]$ can be decomposed as the sum of a regular part $\mu_k^{\rm ac}[u]$ and a 
singular part $\mu_k^s[u]$. As an application of the Theorem 1.1, we prove the following 
theorem. 
 \bigskip
 
\noi
{\bf Theorem 3.2.}\, {\it Let $\Om\subset\Bbb R^n$ be an open set and 
$ u\in\Phi^k(\Om)$, $k>n/2$. Then the absolute continuous 
part of $\mu_k[u]$ is represented by the $k$-Hessian operator $F_k[u]$. That is  
  
\begin{equation}
\label{c1}
\mu_k^{\rm ac}[u]=F_k[u]\,dx\,.
\end{equation}
}

{\it Proof.} Let $u$ be a $k$-convex function, $k>n/2$ and $u_{\ep}$ be the 
mollification of $u$. Then by (\ref{a9}) and properties of mollification 
(see, for example Theorem 1, Section 4.2 in \cite{EG}) it 
follows that $u_{\ep}\in \Phi^k(\Om)\cap C^{\infty}(\Om)$. 
Since $u$ is twice differentiable a.e. (by Theorem 1.1) 
and $u\in W^{2,1}_{\rm loc}(\Om)$ (by Theorem 2.5), we conclude that 
$\na^2 u_{\ep}\to \na^2 u $ in $L^1_{\rm loc}$. Let $\mu_k[u_{\ep}]$ and 
$\mu_k[u]$ are the Hessian measures associated to the functions $u_{\ep}$ and $u$ respectively. 
Then the by weak continuity Theorem 2.3 (Theorem 1.1, \cite{TW2}), $\mu_k[u_{\ep}]$ converges to $\mu_k[u]$ in measure 
and $\mu_k[u_{\ep}]=F_k[u_{\ep}]\,dx$. It follows that for any compact set $E\subset\Om$, 
\begin{equation}
 \label{c3}
 \mu_k[u](E)\geq\limsup_{\ep\to 0} \mu_k[u_{\ep}](E)=\limsup_{\ep\to 0}\int_{E}F_k[u_{\ep}]\,.
 \end{equation}
Since $F_k[u_{\ep}]\geq 0$ and $F_k[u_{\ep}](x) \to F_k[u](x)$ a.e., by Fatou's lemma, for 
every relatively compact measurable subset $E$ of $\Om$, we have 
\begin{equation}
\label{c2}
\int_{E}F_k[u]\leq \liminf_{\ep\to 0}\int_{E} F_k[u_{\ep}]\,.
\end{equation}
Therefore by Theorem 3.1, \cite{TW2}, it follows that $F_k[u]\in L^1_{\rm loc}(\Om)$. Let 
$\mu_k[u]=\mu_k^{\rm ac}[u]+ \mu_k^{\rm s}[u]$, where $\mu_k^{\rm ac}[u]=h\,dx$, 
$h\in L^1_{\rm loc}(\Om)$ and $\mu_k^{\rm s}[u]$ is the singular part supported on 
a set of Lebesgue measure zero. We would like to prove that $h(x)=F_k[u](x)$ ${\cal L}^n$ 
a.e. $x$. By taking 
$E:=\overline B(x,r)$, from (\ref{c3}) and (\ref{c2}), we obtain 
\begin{equation}
\label{c4}
\aint_{\overline B(x,r)}F_k[u]\,dy\leq\frac{\mu_k[u](\overline B(x,r))}
{{\cal L}^n(B(x,r))}=\aint_{\overline B(x,r)}h\,dy +\frac{\mu_k^{\rm s}[u](\overline B(x,r))}
{{\cal L}^n(B(x,r))}\,.
\end{equation}
Hence by letting $\ep\to 0$, we obtain 
\begin{equation}
\label{c5}
F_k[u](x)\leq h(x) \,\,\,\,\,\,\,{\cal L}^n\,\,\,{\rm a.e.}\,\,x. 
\end{equation}
To prove the reverse inequality, let us recall that $h$ is the density of the absolute 
continuous part of the measure $\mu_k[u]$, that is for ${\cal L}^n\,\,\,{\rm a.e.}\,\,x$ 
\begin{equation}
\label{c6}
h(x)=\lim_{r\to 0}\frac{\mu_k^{\rm ac}[u](\overline B(x,r))}
{{\cal L}^n(B(x,r))}=\lim_{r\to 0}\frac{\mu_k[u](\overline B(x,r))}
{{\cal L}^n(B(x,r))}\,.
\end{equation}
Since $\mu^s_k[u]$ is supported on a set of Lebesgue measure zero, 
$$
\mu^s_k[u](\pa B(x,r))=0,\,\,\,\,\,\, {\cal L}^1\,\,\, 
{\rm a.e.}\,\, r>0.
$$
Therefore by the weak continuity of $\mu_k[u_\ep]$ 
(see, for example Theorem 1, Section 1.9 \cite{EG}), we conclude that 
\begin{equation}
\label{c7}
\lim_{\ep\to 0}\mu_k[u_\ep](B(x,r))=\mu_k[u](B(x,r)),\,\,\,\,\,\, {\cal L}^1\,\,\, 
{\rm a.e.}\,\, r>0.
\end{equation}  
Let $\de>0$, then for $\ep<\ep^{\prime}=\ep(\de)$ and for ${\cal L}^1$  
a.e. $r>0$, ${\cal L}^n$ a.e. $x$ 
 \begin{align}
\label{c8}
h(x)&\leq \lim_{r\to 0}\frac{(1+\de)\mu_k[u_\ep](B(x,r))}
{{\cal L}^n(B(x,r))}\nonumber\\
&=(1+\de)\lim_{r\to 0}\aint_{B(x,r)}F_k[u_\ep]\,dy\nonumber\\
&=(1+\de)F_k[u_\ep](x)
\end{align}
By letting $\ep\to 0$ and finally $\de\to 0$, we obtain
$$
h(x)\leq F_k[u](x), {\cal L}^n\,\,\,{\rm a.e.}\,\,x.
$$
This proves the theorem. \qed
 
 \section*{Acknowledgement}
 We would to thank Xu-Jia Wang for many stimulating discussions during this work.

\end{document}